\documentclass[12pt]{article}
\usepackage{amssymb, amsmath}

\newtheorem{theorem}{Theorem}[section]

\newtheorem{definition}{Definition}[section]

\newcommand{\vare}{\varepsilon}
\newcommand{\n}{\nonumber}

\newcommand{\xb}{\bar{x}}
\newcommand{\xbb}{|\bar{x}|}
\newcommand{\si}{\sigma_R (|x|)}

\newcommand{\s}{\sigma}

\newcommand{\bb}{\begin{equation}}
\newcommand{\ee}{\end{equation}}
\newcommand{\bq}{\begin{eqnarray}}
\newcommand{\eq}{\end{eqnarray}}
\newcommand{\bqn}{\begin{eqnarray*}}
\newcommand{\eqn}{\end{eqnarray*}}

\begin{document}
\title{Conditions on the pressure for vanishing velocity in the incompressible fluid flows
in $\Bbb R^N$}
\author{ Dongho Chae
\\
 Department of Mathematics\\
  Sungkyunkwan
University\\
 Suwon 440-746, Korea\\
 e-mail: {\it chae@skku.edu }}
 \date{}
\maketitle
\begin{abstract}
In this paper we derive various sufficient conditions on the
pressure for vanishing velocity in the incompressible Navier-Stokes
and the
Euler equations in $\Bbb R^N$.
{\bf keywords:} Euler equations, Navier-Stokes equations, vanishing
conditions
\end{abstract}

\section{Introduction}
 \setcounter{equation}{0}
 \subsection{The Navier-Stokes and the Euler equations in $\Bbb
 R^N$}
We are concerned  on  the Navier-Stokes and the Euler equations for
incompressible fluid flows on $\Bbb R^N$, $N\in \Bbb N, N\geq 2$.
\[
\mathrm{ (NS, E)}
 \left\{ \aligned
 &\frac{\partial v}{\partial t} +(v\cdot \nabla )v =-\nabla p +\nu \Delta
 v
 \quad (x,t)\in \Bbb R^N\times (0, \infty) \\
 & \textrm{div }\, v =0 , \quad (x,t)\in\Bbb R^N\times (0,
 \infty)\\
  &v(x,0)=v_0 (x), \quad x\in \Bbb R^N
  \endaligned
  \right.
  \]
  where $v(x,t)=(v^1 (x,t), \cdots, v^N (x,t))$ is the velocity, $p=p(x,t)$ is the
  pressure, $\nu \geq 0$ is the viscosity. The Navier-Stokes
  system(NS) corresponds to $\nu >0$, while the Euler system(E)
  corresponds to $\nu =0$.
Given $a,b \in \Bbb R^N$, we denote by $a\otimes b $  the $N\times
N$ matrix with $(a\otimes b)_{ij}=a_ib_j$.
  For two $N\times N$ matrices $A$ and $B$ we denote $A:B=\sum_{i,j=1}^N A_{ij} B_{ij}$.
 Given $m\in \Bbb N\cup \{0\}, q\in [1, \infty]$, we introduce
 $$W^{m,q}_\s (\Bbb R^N):= \left\{ v\in [W^{m,q} (\Bbb R^N)]^N,\,\,\mathrm{ div}\, v=0 \right\},
 $$
 where $W^{m,q} (\Bbb R^N)$ is the standard Sobolev space on $\Bbb R^N,$ and
 the derivatives in  the operation of div $(\cdot)$ are in the sense of
 distribution. In particular, $H^m_\s (\Bbb R^N):=W^{m,2}_\s (\Bbb
\Bbb R^N)$ and $L^q_\s (\Bbb R^N):=W^{0, q}_\s (\Bbb R^N)$.
 Similarly, given $q\in [1, \infty]$, we use $L^q _{ loc, \s } (\Bbb R^N)$ to denote the class of solenoidal
 vector fields, which belongs to $[L^q _{loc} (\Bbb R^N)]^N$.
  In $\Bbb R^N$ we define  weak
 solutions of the Navier-Stokes(Euler)
equations as follows.
  \begin{definition} We say that a pair $(v,p)\in L^2 (0, T; L^2_{loc,\s} (\Bbb R^N))\times L^1
  (0, T; L^1_{loc}(\Bbb R^N ))$ is a weak solution of $(NS,E)$
  on $ \Bbb R^N\times(0, T )$
  if
  \bq\label{11}
  \lefteqn{-\int_{0}^{T} \int_{\Bbb R^N} v(x,t) \cdot \phi(x)\xi ' (t)dxdt -\int_{0}^{T}
  \int_{\Bbb R^N} v(x,t)
  \otimes v(x,t):\nabla \phi (x) \xi (t) dxdt}\hspace{.in} \n \\
  && = \int_{0} ^{T}\int_{\Bbb R^N} p(x,t)\mathrm{ div }\, \phi (x)
  \xi (t)dxdt +\nu \int_{0}^{T} \int_{\Bbb R^N} v(x,t)\cdot \Delta \phi (x) \xi (t)dxdt,\n \\
  \eq
    and
   \bb\label{11a}
  \int_0 ^T \int_{\Bbb R^N} v(x,t)\cdot \nabla \psi (x)\xi(t) dxdt=0
  \ee
   for all $\xi \in C^1 _0 (0, T)$, $\phi \in [C_0 ^\infty (\Bbb R^N
   )]^N$ and $\psi \in C^\infty_0 (\Bbb R^N)$.
  \end{definition}
  Notice that in our definition of weak solution of the Navier-Stokes
  equations  the condition for the velocity is weaker than that for the standard definition, introduced by
  Leray(\cite{ler}), which require $v\in L^\infty (0, T; L_{\s} ^2
  (\Bbb R^N))\cap L^2(0, T; H^1_\s (\Bbb R^N))$.
 If we choose, in particular,  $\phi=\nabla h$ in (\ref{11}), using (\ref{11a}), we
have
 $$
  \int_0 ^T \int_{\Bbb R^N} \xi(t)p(x,t) \Delta h(x)dxdt=-\sum_{j,k=1}^N \int_0 ^T
  \int_{\Bbb R^N} \xi(t)v_j(x,t)v_k (x,t)\partial_j\partial_k h(x)
  dxdt
  $$
  for all $h\in C_0 ^\infty (\Bbb R^N)$ and $\xi \in C^1 _0 (0, T)$,
  from which we deduce
\bb\label{11b}
  \int_{\Bbb R^N} p(x,t) \Delta h(x)dx=-\sum_{j,k=1}^N
  \int_{\Bbb R^N} v_j(x,t)v_k (x,t)\partial_j\partial_k h(x)
  dx
  \ee
for all $h\in C_0 ^\infty (\Bbb R^N)$ and for almost every $t\in [0,
T]$, which is the weak formulation of the well-known relation
between the pressure and the velocity,
 \bb\label{11c} \Delta
p=-\sum_{j,k=1}^N\partial_j\partial_k (v_j v_k)
 \ee
 for the solution of (NS,E).
Our purpose in this paper is to derive sufficient conditions of  the
pressure leading to vanishing velocity from (\ref{11b}). We note
that our study is different to the Liouville type properties of the
weak solution of the Navier-Stokes equations, as investigated in
\cite{koc}, which derive triviality of solution starting from
certain decay assumptions on the solution. Also, we seek the
pressure properties and conditions leading to the triviality of
solutions, which is different from the study of regularity
conditions of the pressure, as done in \cite{ser, ber, cha1, str}.
Our first theorem is the following.
 \begin{theorem}
Let $(v,p)$ be a weak solution of the incompressible
Navier-Stokes(Euler) equations. Let $t\in \{ s\in  [0,\infty) \,| \,
|v(\cdot,s)|^2+ |p(\cdot, s)|\in L^1 (\Bbb R^N)\}$. Then,
 \bb\label{06}
 \int_{\Bbb
 R^{N-1} }p (x,t)d\mathbf{x}'_{j} =-\int_{\Bbb
 R^{N-1} }|v_j(x,t)|^2d\mathbf{x}'_{j}\leq 0
 \ee
for almost every $x_j \in \Bbb R$ and for all $j=1,\cdots, N$, where
we denoted $\mathbf{x}'_{j}:=(x_1, ..., x_{j-1}, x_j, ..., x_N)$.
 \end{theorem}
 {\em Remark 1.1 } Let $\Bbb S^k:=\{ x\in \Bbb R^{k+1}\, |\, |x|=1\}$. Given  $(\xi, x_0) \in \Bbb S^{N-1}\times \Bbb
 R^N$, we define
the  hyperplane $\Pi(\xi, x_0) =\{ x\in \Bbb R^N\, |\,
 \xi \cdot (x-x_0)=0\}$. Then, the formula (\ref{06}) can be generalized
 straightforwardly as follows
 \bb
 \int_{
\Pi(\xi, x_0)}\,p \,dx =-\int_{ \Pi(\xi, x_0)}|v\cdot \xi |^2dx\leq
0\qquad \forall (\xi, x_0)\in \Bbb
 S^{N-1}\times \Bbb R^N ,
 \ee
 which implies that
 $$\{ x\in \Bbb R^N \, |\, p(x,t)\leq 0\} \cap  \Pi(\xi, x_0)\neq \emptyset \qquad \forall (\xi, x_0)\in \Bbb
 S^{N-1}\times \Bbb R^N$$
 for almost every $t\in [0, T)$.\\

\noindent{\em Remark 1.2 }   Previously Brandolese\cite{bra}(see
also \cite{cha} for an independent result) derived
$$ \int_{\Bbb R^N} p(x,t)dx= -\int_{\Bbb R^N} |v_j(x,t)|^2 dx
$$
for all $j=1,\cdots, N,$ which can be derived
from (\ref{06}) by integration over $\Bbb R$  with respect to $x_j$.\\
\begin{theorem} Let $q\in [2, \frac{2N}{N-1} )$ be given. Suppose $v$ is
a weak solution of $(NS,E)$ on $(0, T)$ such that $v(t)\in C(\Bbb
R^N) \cap L^q_\s (\Bbb R^N)$ and $p(t)\in C(\Bbb R^N) \cap L^{1}
(\Bbb R^N)$ for some $t\in (0, T)$. Then, the following formula
holds for all $R\geq 0$.
   \bq\label{th1}
\lefteqn{(N-1)\int_{\{|x|>R\}}\frac{p(x,t)}{|x|}dx+\int_{\{|x|=R\}}
p(x,t)\,d\sigma} \n
\\
&&=-\int_{\{|x|=R\}} |v^r(x,t)|^2\,d\sigma-\int_R^\infty
\int_{\{|x|=r\}}\frac{|v^{\tau}(x,t)|^2}{|x|} d\sigma dr,
 \eq
 where $v^r$ and $v^{\tau}$ are the normal and the tangential components
 of $v$ on the spherical surfaces respectively.
\end{theorem}
Let $\mathcal{S}\subset [0, T]$ be the set of singularity of the
Leray-Hopf weak solution. It is well-known that $\mathcal{S}$ is the
set of Hausdorff
 dimension less than or equal to $1/2$(\cite{ler}).
Since the solution
 $v(x,t)$ can be identical to a smooth function for every
 $t\in [0, T]\setminus \mathcal{S}$, we may assume without loss of generality that
 our solution $v(x,t)$ below is spatially smooth for our choice of
 $t\in[0, T]\setminus
 \mathcal{S}$, which is called the time of regularity.
\begin{theorem}
Let
  $v\in L^\infty(0, T; L^2 _\s(\Bbb R^N))\cap L^2 (0, T; H^1_\s (\Bbb R^N))$ be a
 Leray-Hopf weak solution of $(NS)$, and $t\in \{ s\in [0, T])\setminus \mathcal{S}\, |\,
 p(\cdot, s) \in  L^{1} (\Bbb R^N)\}$.
 Then,
 \bb\label{th13a}
 \int_{\{|x|\geq  R\}} \frac{p(x,t)}{|x|}dx\leq  0\qquad\forall R\geq
 0.
 \ee
Moreover, if there exist $R\geq 0$
 and $t\in (0, T)\setminus \mathcal{S}$
 such that
\bb\label{th13b}
 \int_{\{|x|\geq  R\}} \frac{p(x,t)}{|x|}dx= 0,
 \ee
then,  $v(\cdot,t)=0$ for almost every $t\in [0, T)$.
\end{theorem}

\subsection{Axisymmetric flows in $\Bbb
R^N$}

 Let $N\geq 3$. We consider the system (NS, E) with
the axial symmetry in $\Bbb R^N$. Let us denote $x=(\bar{x}, x_N)\in
\Bbb R^{N-1}\times \Bbb
 R$. We use spherical coordinates to represent $\Bbb R^{N-1}$,
 namely
 $$\bar{x}=(\rho, \theta_1, \cdots, \theta_{N-2})\in \Bbb R_+\times
 \Bbb S^{N-2},
 $$
 where $\rho=|\bar{x}|$, and $\theta_1, \cdots, \theta_{N-2}$ are
 the angular variables on $\Bbb S^{N-2}$. Let
 $\{ e_\rho,
 e_{\theta_1}, \cdots, e_{\theta_{N-2}},  e_N\}
 $
 be the canonical basis of this generalized cylindrical coordinate
 system.
 In this coordinate system the velocity field can be
 represented as
 $$
 v=v^\rho e_\rho+v^{\theta_1}
 e_{\theta_1}+\cdots+ v^{\theta_{N-2}}e_{\theta_{N-2}}+v_N
 e_N.
 $$
 By an axisymmetric flow in $\Bbb R^N$ we mean that each of the
 components $v^\rho , \cdots, v^N$ does not depend on the angular
 variables $\theta_1, \cdots, \theta_{N-2}$.
  In this case the swirl free flow means that $v^{\theta_1}=\cdots
  =v^{\theta_{N-2}}=0$.
 \begin{theorem} Let $N\geq 3$ and $(v,p)$ be an axially symmetric
   classical solution of $(NS,E)$ without swirl
  on $ \Bbb R^N \times(0, T )$.  We define
 $$
I(\rho,  t):=\int_{\Bbb R} \left\{(v^\rho(x,t)) ^2 +p(x,t)\right\}
dx_N.
 $$
  Then, $I(\rho, t)$ decays monotonically in radial direction. More precisely, we have
\bb\label{th21}
  I(\rho_2,t)-I(\rho_1,t) =-(N-2)
 \int_{\rho_1}
 ^{\rho_2}\int_{\Bbb R} \frac{(v^\rho) ^2 }{\rho}d x_N d\rho\leq 0
\ee
 for all $\rho_2>\rho_1 >0$.  In particular, if $I(\cdot ,t)$  decays to zero at infinity,
 then we have
 \bb\label{th22}
\int_{\Bbb R} p(0, x_N,t)dx_N=(N-2)\int_{0 }^\infty
 \int_{\Bbb R}\frac{(v^\rho) ^2 }{\rho}dx_Nd\rho.
 \ee
 Thus, $\int_{\Bbb R} p(0, x_N,t)dx_N=0$ implies $v=0$ in this case.
\end{theorem}
\noindent{\em Remark 1.3 } Similar results hold if we have periodic
condition in the direction of $x_N$ with period $L$. In this case
(\ref{th22}) reduces to
  \bb\label{th22a} \int_{0} ^L p(0,
x_N,t)dx_N=(N-2)\int_0 ^L \int_{\Bbb R}
 \frac{(v^\rho) ^2 }{\rho}dx_Nd\rho.
 \ee

\section{Proof of the Main Results}
 \setcounter{equation}{0}

\noindent{\bf Proof of Theorem 1.1 }
 Let us consider a radial cut-off function $\sigma\in C_0
^\infty(\Bbb R^N)$ such that
 \bb\label{07}
   \sigma (x)=\sigma(|x|)=\left\{ \aligned
                  &1 \quad\mbox{if $|x|<1$}\\
                     &0 \quad\mbox{if $|x|>2$},
                      \endaligned \right.
 \ee
and $0\leq \sigma  (x)\leq 1$ for $1<|x|<2$. We set $\sigma_R
(x)=\sigma (\frac{x}{R})$ for $R>0$.
  Given $ m\in
\{ 1, \cdots , N\}$, we set $ w(x)= e^{i\xi_m x_m}$, and $h (x)=
w(x)\sigma_R (x)$ in (\ref{11b}), and pass $R\to \infty$.
 Then,  from the hypothesis $|v(\cdot, t)|^2 +|p(\cdot, t)|\in L^1 (\Bbb
 R^N)$ one can apply the dominated convergence theorem to  obtain \bqn
  0& =&  \int_{\Bbb R^N} \left\{\sum_{j,k=1}^Nv_j (x,t) v_k
  (x,t)\partial_j\partial_k w(x)
  +p(x,t)\Delta w(x)\right\}
  dx\n \\
  &=&-\xi_m ^2 \int_{\Bbb R^N} (|v_m (x,t)|^2 +p(x,t)) e^{i\xi_m
 x_m
 }dx\n \\
 & =& -\xi_m ^2 \int_{-\infty} ^{+\infty} \left\{ \int_{\Bbb R^{N-1}}
(|v_m (x,t)|^2 +p(x,t)) d\mathbf{x}'_m\right\}e^{i\xi_m x_m
 }dx_m,
 \eqn
 which shows that
 $
  \hat{f}(\xi_m)=0
 $
for all $\xi_m \neq 0$, where we set
 $
  f(x_m )=\int_{\Bbb R^{N-1}}
(|v_m (x,t)|^2 +p(x,t))d\mathbf{x}'_m.
 $
 Since $ \hat{f}\in C_0 (\Bbb R)$, we find that
 $\hat{f}(\xi_m )=0 $  for all $\xi_m \in \Bbb R$ by continuity.
   Hence $ f(x_m)=0$ for almost every $ x_m \in \Bbb R$. $\square$\\
\ \\

\noindent{\bf Proof of Theorem 1.2 } Let us consider the function
$\eta \in C^\infty _0 ( \Bbb R^N)$ such that
$$\eta(x)=\left\{ \aligned &c\exp \left(\frac{1}{|x|^2 -1}\right), &
|x|<1,\\
 &0, & |x|>1,
 \endaligned \right.
 $$
 where the constant $c$ is the normalizing constant so that $\int_{\Bbb R} \eta
 (x)dx=1$. For $\vare >0$ we define the sequence of mollifiers
 $\eta_\vare (x)=\frac{1}{\vare} \eta (\frac{x}{\vare})$.
Given $R_2 > R_1>0$ we choose $\vare$ below so that $0<\vare<\min\{
R_1, \frac{R_2-R_1}{2}\}$. For such $R_1, R_2,\vare, \eta_\vare
(\cdot )$ we define $\varphi_{R_1, R_2, \vare}(r)=\varphi_{R_1, R_2,
\vare}(|x|)$ by
 \bb\label{13}
  \varphi_{R_1,R_2,\vare}(r)= \int_0 ^{r}\int_0 ^\sigma \left\{\eta_\vare (s-R_1)-\eta_\vare (s-R_2
  )\right\} dsd\sigma.
 \ee
 We observe that $\varphi_{R_1,
R_2,\vare} (r)$ satisfies
 \bb\label{14}
  \varphi'_{R_1,R_2,\vare}(r) \to \chi_{\{ R_1< r<R_2\}}
  (r),\quad\quad
    \varphi^{\prime\prime}_{R_1,R_2,\vare}(r) \to \delta(r-R_1)-\delta (r-R_2)
  \ee
 as $\vare\to 0$   in the sense of distribution. Moreover, for $g\in
C_0 (\Bbb R)$ we have
 \bb\label{14a}
 \lim_{\vare\to 0}\int_0 ^\infty \varphi'_{R_1,R_2,\vare}(r)g(r)dr
 =\int_{R_1} ^{R_2}g(r)dr
 \ee
 for all $R_2>R_1>0$, and
 \bb\label{14b}
 \lim_{\vare\to 0}\int_0 ^\infty
 \varphi^{\prime\prime}_{R_1,R_2,\vare}(r)g(r)dr
 =g(R_1)-g(R_2).
 \ee
Let us choose the test function
 $h$ in (\ref{11b}) as
 \bb\label{15}
  h(x)=  \varphi_{R_1,R_2, \vare} (|x|)\si,
 \ee
 where $\si$ is the smooth cut-off function introduced in the proof
 of Theorem 1.1. Then, after passing $R\to\infty$, and using the
 dominated convergence theorem, (\ref{11b}) becomes
  \bq\label{15a}
 \lefteqn{ \int_{\Bbb R^N}   v_j (x,t) v_k (x,t)\partial_j\partial_k \{ \varphi_{R_1,R_2, \vare}
 (|x|)\} dx} \n \\
 &&=-\int_{\Bbb R^N}  p(x,t) \Delta\{ \varphi_{R_1,R_2, \vare}
 (|x|)\} dx
\eq
 for every $t\in [0, T]$.
 Then, since
 $$ \partial_j\partial_k h(|x|)= \left(\frac{\delta_{jk}}{|x|}
 -\frac{x_jx_k}{|x|^3} \right) h' (|x|) +\frac{x_j x_k} {|x|^2}
 h^{\prime\prime} (|x|),
 $$
 we obtain from (\ref{15a}) that
 \bq\label{16}
\lefteqn{0=\int_{\Bbb R^N} \varphi_{R_1,R_2, \vare}^{\prime\prime}
(|x|) \left\{\sum_{j,k=1}^N v_j v_k \frac{x_jx_k}{|x|^2}
+p(x,t)\right\} \,dx
}\hspace{.0in}\n\\
&& + \int_{\Bbb R^N}\varphi_{R_1,R_2, \vare}^{\prime} (|x|)
\left\{\frac{|v|^2}{|x|} -\sum_{j,k=1}^N v_j v_k
\frac{x_jx_k}{|x|^3}\right.\n \\
&&\hspace{2.in} \left.+(N-1)\frac{p(x,t)}{|x|} \right\}  \,dx,
 \eq
 which can be rewritten as
 \bq\label{17a}
 \lefteqn{ \int_0^\infty \varphi_{R_1,R_2, \vare}^{\prime\prime}
(r) \int_{\{|x|=r\}}\left\{ (v^r(x,t))^2  +p(x,t)\right\} \,d\sigma dr}\n \\
&&=-\int_0^\infty \varphi_{R_1,R_2, \vare}^{\prime} (r)
\int_{\{|x|=r\}} \left\{\frac{(v^{\tau}(x,t))^2}{r}+(N-1)\frac{p(x,t)}{r} \right\}\,d\sigma  dr\n \\
 \eq
Passing $\vare \to 0$ in (\ref{17a}), and using the facts
(\ref{14a}) and (\ref{14b}), we deduce
 \bq\label{17b}
 \lefteqn{
\int_{\{|x|=R_1\}}( (v^r(x,t))^2  +p(x,t)) \,d\sigma-\int_{\{|x|=R_2\}}( (v^r(x,t))^2  +p(x,t)) \,d\sigma }\n \\
&&=-\int_{R_1}^{R_2} \int_{\{|x|=r\}}
\left\{\frac{(v^{\tau}(x,t))^2}{|x|}+(N-1)\frac{p(x,t)}{|x|}
\right\}\,d\sigma dr.
 \eq
 Now we observe
 \bq\label{19}
 \lefteqn{\left|\int_{\{R_1<|x|<R_2\}} \frac{|v^{\tau} (x,t)|^2}{|x|}
dx\right|\leq C \int_{\{|x|>R_1\}} \frac{|v(x,t)|^2}{|x|}dx }\n \\
&&\quad\leq C\left(\int_{\Bbb R^N}|v(x,t)|^{q}
dx\right)^{\frac{2}{q}} \left(\int_{\{|x|>R_1\}}
\frac{1}{|x|^{\frac{q}{q-2}}} dx\right)^{\frac{q-2}{q}}\n \\
&&\quad \leq C  \|v(t)\|_{L^q}^2 \left(\int_{R_1} ^\infty r
^{N-1-\frac{q}{q-2}} dr\right)^{\frac{q-2}{q}} <\infty
 \eq
for $q\in [2, \frac{2N}{N-1} )$, and in the limiting case $q=2$ we
have in mind the obvious estimate,
$$ \int_{\{|x|>R_1\}} \frac{|v(x,t)|^2}{|x|}dx\leq \frac{1}{R_1}
\|v(t)\|_{L^2}^2.
$$
On the other hand,
  \bq\label{20}
\lefteqn{\left|\int_{\{R_1<|x|<R_2\}} \frac{p(x,t)}{|x|}dx \right|
\leq
\int_{\{|x|>R_1\}} \frac{|p(x,t)|}{|x|}dx}\n \\
&&\quad\leq C\|p(t)\|_{L^{\frac{q}{2}}} \left(\int_{\{|x|>R_1\}}
\frac{1}{|x|^{\frac{q}{q-2}}} dx\right)^{\frac{q-2}{q}} \n \\
&&\quad\leq C \|v(t)\|_{L^q}\left( \int_{R_1} ^\infty r
^{N-1-\frac{q}{q-2}} dr\right)^{\frac{q-2}{q}} <\infty
 \eq
 for $1<q< \frac{2N}{N-1}$,  where we used the
Calderon-Zygmund inequality for the relation $p=-\sum_{j,k=1}^N
\partial_j\partial_k \Delta^{-1} (v_j v_k )$.
 Passing $R_2 \to \infty$ in (\ref{17b}), using the
dominated convergence theorem, which is justified by the facts
(\ref{19}) and (\ref{20}), we obtain (\ref{th1}).
 $\square$\\
 \ \\
 \noindent{\bf Proof of Theorem 1.3 } We first observe  the
 inclusion relation
 $$L^\infty(0, T;  L^2(\Bbb R^N))\cap L^2 (0, T; H^1 (\Bbb R^N))\subset L^2 (0, T; L^q (\Bbb R^N)),
 $$
 for any $q\in [2, \frac{2N}{N-2}]$ when $N>2$, and for all $q\in
 [2, \infty)$ when $N=2$.
Thus, for the Leray-Hopf weak solution $v$ we have
 $$v\in L^2 (0, T; L^q (\Bbb R^N)), \quad p\in L^1 (0, T; L^\frac{q}{2} (\Bbb
 R^N))$$
 for $2< q< \frac{2N}{N-1}$, $T>0$.
 Let $t_1\in [0, \infty)\setminus \mathcal{S}$ be fixed. Then,
 $v(\cdot, t_1)$ and $p(\cdot, t_1)$ are smooth functions. Now suppose
 there exists $R\geq 0$ such that
\bb\label{de}
 \int_{\{|x|\geq  R\}} \frac{p(x,t_1)}{|x|}dx>0 .
 \ee
Then, the formula (\ref{th1}) implies that
$$\int_{\{|x|=R\}}
p(x,t_1)d\sigma < 0. $$
  Let us define
$$R_1=\inf\left\{r> R\, \Big|\, \int_{\{|x|=r\}} p(x,t_1)d\sigma > 0\right\}.
$$
By the hypothesis (\ref{de}) and the identity
$$\int_{\{|x|\geq  R\}} \frac{p(x,t_1)}{|x|}dx =\int_R ^\infty \frac{1}{r}\int_{\{|x|=r\}}
p(x,t_1)d\sigma dr
$$
we have $R< R_1 <\infty$. Moreover,
  \bb\label{pr1}
  \int_{\{|x|=R_1 \}} \frac{p(x,t_1)}{R_1}d\sigma =0,
 \ee
 and, since
 $$\int_{\{ R<|x|<R_1\}} \frac{p(x,t_1)}
 {|x|} dx=\int_{R} ^{R_1}\int_{\{|x|=r\}}\frac{p(x,t_1)}
 {|x|}d\sigma dr \leq 0,$$
 we have
 \bb\label{pr2}
 \int_{\{|x|\geq  R_1\}} \frac{p(x,t_1)}{|x|}dx\geq \int_{\{|x|\geq  R\}}
 \frac{p(x,t_1)}{|x|}dx >0.
 \ee
 Combining (\ref{pr1}) and (\ref{pr2}) with the formula (\ref{th1})
 we have the contradiction. Therefore we obtain (\ref{th13a}).
Next, we suppose there exists $t_2 \in [0, \infty)\setminus
\mathcal{S}$ such that
  \bb\label{pr3}
 \int_{\{|x|\geq  R\}} \frac{p(x,t_2)}{|x|}dx=0
 \ee
 Then,
 by the formula (\ref{th1}) we have
 \bb\label{pr4}
\int_{\{|x|= R\}} \frac{p(x,t_2)}{|x|}d\sigma \leq 0
 \ee
We claim that there exists $R_2 \in [R, \infty)$ such that
  \bb\label{21}
  \int_{\{|x|=R_2\}} |v^r (x,t_2)|^2d\sigma (x)+\int_{R_2}^\infty\int_{\{|x|=r\}}
  \frac{|v^{\tau} (x,t_2)|^2}{|x|} d\sigma dr = 0.
  \ee
 In the case of equality in (\ref{pr4}), combining this with
 (\ref{pr3}) we have (\ref{21}) by the formula (\ref{th1}).
 If the strict inequality holds in (\ref{pr4}), then, we define
  $R_2\in [R, \infty)$ by
$$ R_2=\inf\left\{r> R\, \Big|\, \int_{\{|x|=r\}} p(x,t_2)d\sigma > 0\right\}.
$$
Note that (\ref{pr3}) implies $R<R_2 <\infty$. Then, by (\ref{pr3})
and continuity of $p(\cdot, t_2)$ we have
  \bb\label{pr5}
 \int_{\{|x|\geq  R_2 \}} \frac{p(x,t_2)}{|x|}dx>0,
 \ee
which contradicts (\ref{th13a}). Thus, we have established
(\ref{21}), from which we immediately have
 \bb\label{22} v^r (x,t_2)=0 \quad  \mbox{on $\{|x|=R_2\}$},
 \ee
   and
 \bb\label{23} v(x,t_2)=v^r (x,t_2)\quad \mbox{on $\{|x|>R_2\}$}.
  \ee
 In particular, (\ref{23}) implies that
 $$0=\mathrm{div}\, v=\frac{1}{r^{N-1}}\partial_r\{ r^{N-1}
 v^r(x,t_2)\}=0 \quad \mbox{on}\quad \{|x|>R_2\}
 $$
 and hence
 \bb\label{24}
 v^r (x,t_2) =\frac{C}{r^{N-1}}\quad \quad \mbox{on $\{|x|>R_2\}$},
\ee
  where $C=C(\theta, t_2)$ depends only on the angular variables
$\theta \in \Bbb S^{N-1}$ and $t_2$. Passing $r\downarrow R_2$ in
(\ref{24}), and applying the condition (\ref{22}), we obtain $C=0$.
Thus we have $v(\cdot,t_2)=0$ on $\{|x|>R_2\}$. Since $v(\cdot,t_2)$
is a Leray-Hopf weak solution with a compact support on $\Bbb R^N$,
we conclude $v=0$ on $\Bbb  R^N\times[t_2, t_2+\vare) $ for some
$\vare
>0$(see e.g. \cite{sch}). Thus,
we conclude that $v=0$ almost everywhere in $\Bbb R^N\times [0, T)$. $\square$\\
\ \\

 \noindent{\bf Proof of Theorem 1.4 }
We consider
 \bb\label{13}
  \varphi_{R,\vare}(r)=
  \int_0 ^{r}\eta_\vare (s-R) ds,
 \ee
 where $\eta_\vare (\cdot)$ is the approximation of  identity defined the proof
 of
 Theorem 1.2.
Below we denote $\bar{v}:=(v_1, \cdots, v_k, 0)$, which is the
projection
 of $v(x,t)\in \Bbb R^N$ onto $\Bbb R^{N-1}$, and $v^\rho :=v\cdot
 \frac{\xb}{\xbb }$.
 Then, similarly to the proof of Theorem 1.2 we compute
  \bq\label{26a}
  \lefteqn{ \int_{\Bbb R^N}     \varphi_{R,\vare}^{\prime\prime}
  (\xbb)\left\{ \sum_{j,k=1}
  ^{N-1} v_j (x,t)v_k (x,t)\frac{x_j x_k}{|\xb|^2} +
  p(x,t) \right\}\sigma_{R_1} (x_N) dx }\n \\
  &&+\int_{\Bbb R^N}  \varphi_{R,\vare}^{\prime}(\xbb)
  \left\{ \frac{|\bar{v}(x,t)|^2}{\xbb} -\sum_{j,k=1}
  ^{N-1} v_j (x,t)v_k (x,t)\frac{x_j x_k}{|\xb|^3}\right\}\sigma_{R_1}
  (x_N)\,dx\n \\
  &&\quad+(N-2) \int_{\Bbb R^N}
    \varphi_{R,\vare}^{\prime}(\xbb)\frac{p(x,t)}{\xbb}
  \sigma_{R_1} (x_N)\,dx\n \\
  &&\quad= -\frac{1}{R^2_1} \int_{\Bbb R^N} \{v_N (x,t)^2+p(x,t)\}
  \varphi_{R,\vare}(\xbb)  \sigma^{\prime\prime}
(\frac{x_N}{R}) dx\n \\
  &&\quad\quad-\frac{2}{R_1}
  \sum_{j=1} ^{N-1}\int_{\Bbb R^N}v_j  (x,t)
  \frac{x_j}{\xbb}
v_N (x,t)
  \varphi_{R,\vare}(\xbb)\sigma^{\prime}
(\frac{x_N}{R}) dx.
  \eq
 Passing
 $R_1\to \infty$ in (\ref{26a}),  using the dominated convergence theorem
 repeatedly, we obtain that
 \bq\label{27}
&& \int_{\Bbb R^N} \left[  \varphi_{R,\vare}^{\prime}(\xbb) \left(
 \frac{|\bar{v}|^2}{\xbb} -\frac{(\bar{v}\cdot
 \xb)^2}{\xbb^3}\right) +  \varphi_{R,\vare}^{\prime}(\xbb) \frac{(N-2)p}{\xbb}
 \right.\n \\
 &&\hspace{1.5in}\left.+  \varphi_{R,\vare}
 ^{\prime\prime}(\xbb)
 \left(\frac{(\bar{v}\cdot \xb )^2}{\xbb^2} +p \right)\right]dx =0
 \eq
 Let us denote by $\omega_N$ the volume of the unit ball in $\Bbb
 R^N$.
 Using the cylindrical coordinate system, and integrating by part,
 we have
 \bqn
\lefteqn{\int_{\Bbb R^N}  \varphi_{R,\vare}
 ^{\prime\prime}(\xbb)
 \left(\frac{(\bar{v}\cdot \xb )^2}{\xbb^2} +p \right)dx}\n \\
 &&=
(N-1)\omega_{N-1} \int_{\Bbb R} \int_0 ^\infty  \varphi_{R,\vare}
 ^{\prime\prime}(\rho)((v^\rho) ^2 +p)\rho ^{N-2}d\rho dx_N\n \\
 &&=-(N-1)\omega_{N-1}  \int_{\Bbb R} \int_0 ^\infty  \varphi_{R,\vare}
 ^{\prime}(\rho) \partial_\rho ((v^\rho) ^2 +p)\rho ^{N-2}d\rho dx_N\n \\
 &&\qquad\quad -(N-2)(N-1)\omega_{N-1}  \int_{\Bbb R}
 \int_0 ^\infty  \varphi_{R,\vare}
 ^{\prime}(\rho) ((v^\rho) ^2 +p)\rho ^{N-3}d\rho dx_N
 \eqn
Therefore we can rewrite (\ref{27}) as in
 following.
 \bq\label{28}
 0&=&(N-1)\omega_{N-1} \int_{\Bbb R} \int_0 ^\infty\left[\frac{\varphi_{R, \vare}^{\prime} (\rho)}{\rho} (|\bar{v}|^2-v_\rho^2 )
 +\varphi_{R, \vare}^{\prime}(\rho) \frac{(N-2)p}{\rho}\right]\rho^{N-2} d\rho dx_N \n \\
 &&\qquad -(N-1)\omega_{N-1}  \int_{\Bbb R} \int_0 ^\infty\varphi_{R, \vare}
 ^{\prime}(\rho) \partial_\rho ((v^\rho) ^2 +p)\rho ^{N-2}d\rho dx_N\n \\
 &&\qquad\quad -(N-2)(N-1)\omega_{N-1}  \int_{\Bbb R} \int_0 ^\infty\varphi_{R, \vare}
 ^{\prime}(\rho) ((v^\rho) ^2 +p)\rho ^{N-3}d\rho dx_N\n \\
 &=&(N-1)\omega_{N-1}  \int_{\Bbb R} \int_0 ^\infty\left[\frac{\varphi_{R, \vare}^{\prime}
 (\rho)}{\rho} (|\bar{v}|^2-(N-1) (v^\rho)^2 )\right]\rho ^{N-2}d\rho dx_N\n \\
 &&\quad-(N-1)\omega_{N-1} \int_{\Bbb R} \int_0 ^\infty\varphi_{R, \vare}
 ^{\prime}(\rho) \partial_\rho ((v^\rho) ^2 +p)\rho ^{N-2}d\rho dx_N.
 \eq
Passing $\vare \downarrow 0$ in (\ref{28}) and  dividing by
$(N-1)\omega_{N-1} R^{N-2}$, and setting $R=\rho$, we obtain
  \bq\label{29}
\lefteqn{\partial_\rho \int_{\Bbb R} \{p(\rho, x_N,t)+ (v^\rho(\rho,
x_N,t ))
^2\}d x_N}\hspace{.5in}\n \\
&&= \frac{1}{\rho}\int_{\Bbb R}\left\{ |\bar{v}(\rho,x_N,t)|^2-(N-1)
(v^\rho (\rho,x_N,t))^2\right\}d
 x_N\n \\
 &&=\frac{1}{\rho}\int_{\Bbb R}\left\{
 \sum_{j=1}^{N-2}(v^{\theta_j}(\rho,x_N,t))^2-(N-2) (v^\rho(\rho,x_N,t))^2\right\} d
 x_N.\n \\
 \eq
 In the case of the swirl free flows  we have
 \bb\label{30}
\partial_\rho \int_{\Bbb R} (p+ (v^\rho) ^2 )d x_N=
 -\frac{(N-2)}{\rho}\int_{\Bbb R} (v^\rho)^2  dx_N.
 \ee
  Integrating
(\ref{30}) over $(\rho_1, \rho_2)$ we find that
 \bb\label{31}
 I(\rho_2,t)-I(\rho_1,t) =-(N-2)
 \int_{\rho_1}
 ^{\rho_2}\int_{-\infty} ^\infty \frac{(v^\rho) ^2 }{\rho}d x_N d\rho
\ee
  In the case when $I(+\infty, t)=0$, the equality (\ref{th22})
  follows by taking the limits $\rho_2\to \infty$, and then $\rho_1
  \to 0$ in (\ref{31}), observing $v^\rho(\rho,x_N,t) \to 0$ as $\rho\to 0$
  for all $x_N\in \Bbb R$, $t\in [0, T]$ if $v(x,t)$ is
   a classical solution of (NS,E). $\square$\\
$$\mbox{\bf Acknowledgements} $$
The author would like to thank to Prof. P. Constantin, the
communication with whom inspired the proof of Theorem 1.1. He also
would like to the thank to the anonymous referee for careful reading
and helpful suggestions. This work was supported partially by the
NRF grant. no. 2006-0093854.


\begin{thebibliography}{1}
\bibitem{ber}L. C. Berselli and G. P. Galdi, {\it Regularity criteria
involving the pressure for the weak solutions to the Navier-Stokes
equations,} Proc. Amer. Math. Soc., {\bf 130},  no. 12, (2002), pp.
3585-3595.
\bibitem{bra}L. Brandolese and Y. Meyer, {\it On the instantaneous spreading for
the Navier-Stokes system in the whole space,}  Contr. Optim. Calc.
Var. {\bf 8},  (2002), pp. 273-285.
\bibitem{cha} D. Chae, {\it Liouville type of theorems for
the Euler and the Navier-Stokes equations,} Adv.  Math., {\bf 228},
(2011), pp. 2855-2868.
\bibitem{cha1} D. Chae and J. Lee, {\it Regularity criterion in terms
of pressure for the Navier-Stokes equations}, Nonlinear Anal., {\bf
46}, (2001), pp. 727-735.
 \bibitem{koc}G. Koch, N. Nadirashvili, G. Seregin and V. \v{S}ver\'{a}k,
 {\it Liouville theorems for the Navier-Stokes equations and
 applications,}  Acta Math., {\bf 203},  no. 1, (2009), pp. 83-105.
 \bibitem{ler} J. Leray, {\it Essai sur le mouvement d'un fluide
visqueux emplissant l'espace,} Acta Math., {\bf 63}, (1934), pp.
193-248.
\bibitem{sch}M. Schonbek, {\it Lower bounds of rates of decay for
solutions to the Navier-Stokes equations,} J. Amer. Math. Soc., {\bf
4}, no. 3, (1991), pp. 423-449.
 \bibitem{ser}G. Seregin and V. \v{S}ver\'{a}k, {\it
Navier-Stokes equations with lower bounds on the pressure,} Arch.
Ration. Mech. Anal., {\bf  163},  no. 1, (2002),  pp. 65-86.
\bibitem{str}M. Struwe, {\it On a Serrin-type regularity criterion for the
Navier-Stokes equations in terms of the pressure,}  J. Math. Fluid
Mech., {\bf  9}, no. 2, (2007), pp. 235-242.
\end{thebibliography}
\end{document}